\documentclass[12pt]{amsart}

\newtheorem{THM}{Theorem}[section]
\newtheorem{LMA}[THM]{Lemma}
\newtheorem{PROP}[THM]{Proposition}
\newtheorem{CORO}[THM]{Corollary}

\numberwithin{equation}{section}

\usepackage{amsmath}
\usepackage{amssymb}
\usepackage{euscript}

\newcommand{\showon}{\begin{eqnarray*}}
\newcommand{\showoff}{\end{eqnarray*}}

\newcommand{\none}{\varnothing}
\newcommand{\drop}{\smallsetminus}

\newcommand{\goesto}{\rightarrow}
\newcommand{\one}{\boldsymbol{1}}
\newcommand{\zero}{\boldsymbol{0}}

\newcommand{\A}{\EuScript{A}} \renewcommand{\a}{\mathbf{a}}
\newcommand{\B}{\EuScript{B}}\renewcommand{\b}{\mathbf{b}}
 \newcommand{\CC}{\mathbb{C}}
\renewcommand{\c}{\mathbf{c}}
 \renewcommand{\d}{\mathbf{d}}

\newcommand{\F}{\EuScript{F}} 

\newcommand{\G}{\EuScript{G}}

\newcommand{\J}{\EuScript{J}}

\renewcommand{\L}{\EuScript{L}}
\newcommand{\M}{\EuScript{M}} \newcommand{\m}{\mathbf{m}}
\newcommand{\N}{\EuScript{N}}

\renewcommand{\P}{\EuScript{P}} 

\newcommand{\Q}{\EuScript{Q}} 
\newcommand{\q}{\mathbf{q}}
 \newcommand{\RR}{\mathbb{R}}
 \renewcommand{\S}{\EuScript{S}}

\newcommand{\U}{\EuScript{U}}\renewcommand{\u}{\mathbf{u}}
\renewcommand{\v}{\mathbf{v}}
\newcommand{\w}{\mathbf{w}}

\newcommand{\Y}{\EuScript{Y}} \newcommand{\y}{\mathbf{y}}

\begin{document}

\title[Rayleigh matroids]{Rayleigh Matroids}

\author{YoungBin Choe}
\address{Department of Combinatorics and Optimization\\
University of Waterloo\\
Waterloo, Ontario, Canada\ \ N2L 3G1}
\email{\texttt{ybchoe@math.uwaterloo.ca}}

\author{David G. Wagner}
\address{Department of Combinatorics and Optimization\\
University of Waterloo\\
Waterloo, Ontario, Canada\ \ N2L 3G1}
\email{\texttt{dgwagner@math.uwaterloo.ca}}


\begin{abstract}
Motivated by a property of linear resistive electrical networks,
we introduce the class of Rayleigh matroids.  These form a subclass
of the balanced matroids defined by Feder and Mihail \cite{FM} in 1992.
We prove a variety of results relating Rayleigh matroids to other
well--known classes -- in particular, we show that a binary matroid
is Rayleigh if and only if it does not contain $\S_{8}$ as a minor.
This has the consequence that a binary matroid is balanced if and
only if it is Rayleigh, and provides the first complete proof in
print that $\S_{8}$ is the only minor--minimal binary
non--balanced matroid, as claimed in \cite{FM}.  We also give an
example of a balanced matroid which is not Rayleigh.
\end{abstract}

\maketitle

\section{Introduction.}

(For explanation of any undefined terms, we refer the reader to
Oxley's book \cite{Ox}.)

In 1992, Feder and Mihail \cite{FM} introduced the concept of a
balanced matroid in relation to a conjecture of Mihail
and Vazirani \cite{MV} regarding expansion properties of 
one--skeletons of $\{0,1\}$--polytopes.
(Unfortunately, the term ``balanced'' has also been
used for matroids with at least three other meanings \cite{B,CM,HL}.)
Let $\M$ be a matroid with ground--set $E$.  For disjoint
subsets $I,J$ of $E$, let $\M_{I}^{J}$ denote the minor of $\M$ 
obtained by contracting $I$ and deleting $J$, and let $M_{I}^{J}$ 
denote the the number of bases of $\M_{I}^{J}$.  Feder and
Mihail say that $\M$ is \emph{negatively correlated} provided that
for every $e,f\in E$ with $e$ not a loop,
$$\frac{M_{f}}{M} \geq \frac{M_{ef}}{M_{e}},$$
and that $\M$ is \emph{balanced} provided that every minor of $\M$
is negatively correlated.  Since $M_{e}=M_{ef}+M_{e}^{f}$,
$M_{f}=M_{ef}+M_{f}^{e}$, and $M=M_{ef}+M_{e}^{f}+M_{f}^{e}+
M^{ef}$, the inequality above is equivalent to
$$\Delta M\{e,f\}:=M_{e}^{f}M_{f}^{e}-M_{ef}M^{ef}\geq 0.$$
We briefly review the literature on balanced matroids in Section 2.

Stemming from a collaboration with James Oxley and Alan Sokal
\cite{COSW}, we were motivated to consider the following similar
condition on a matroid $\M$ with ground--set $E$.  Fix indeterminates
$\y:=\{y_{e}:\ e\in E\}$ indexed by $E$, and for disjoint subsets
$I,J\subseteq E$ let $M_{I}^{J}(\y):=\sum_{B}\y^{B}$, with the sum
over all bases $B$ of $\M_{I}^{J}$ and with $\y^{B}:=\prod_{e\in B}y_{e}$.
We say that $\M$ is a \emph{Rayleigh} matroid provided that whenever
$y_{c}>0$ for all $c\in E$, then for every pair of distinct $e,f\in E$,
$$\Delta 
M\{e,f\}(\y):=M_{e}^{f}(\y)M_{f}^{e}(\y)-M_{ef}(\y)M^{ef}(\y)\geq 0.$$
We call the polynomial $\Delta M\{e,f\}(\y)$ the \emph{Rayleigh
difference of $\{e,f\}$ in $\M$}.  This terminology is motivated by
the Rayleigh monotonicity property of linear resistive electrical
networks, as explained in Section 3.
The main results of Section 3 are as follows.\\
$\bullet$\   The class of Rayleigh matroids is closed by taking
duals and minors.\\
$\bullet$\ Every Rayleigh matroid is balanced.\\
$\bullet$\ The class of Rayleigh matroids is closed by taking
$2$--sums.\\
$\bullet$\ The class of balanced matroids is closed by taking
$2$--sums if and only if every balanced matroid is Rayleigh.\\
$\bullet$\ A binary matroid is Rayleigh if and only if it does
not contain $\S_{8}$ as a minor.\\
$\bullet$\ A binary matroid is balanced if and only if it is 
Rayleigh.\\
These results were motivated by similar claims for balanced matroids
for which complete published proofs are not available.

In Section 4 we discuss another class of matroids -- the
``half--plane property'' matroids, or HPP matroids for short.
This class was, in part, the object of study in our collaboration
with Oxley and Sokal \cite{COSW}.  We extend a theorem of Godsil
\cite{Go} (itself a refinement of a theorem of Stanley \cite{St})
from the class of regular matroids to the more general class of HPP
matroids. The following consequence of this is the main result of
Section 4:\\
$\bullet$ Every HPP matroid is a Rayleigh matroid.\\
In proving this we identify a spectrum of conditions between
these two extremes.  

In Section 5 we discuss some more specific examples.
On the positive side:\\ 
$\bullet$\ All sixth--root of unity matroids are HPP matroids.
(This is from \cite{COSW}.)  In particular, all regular matroids
(hence all graphs) are HPP matroids, and hence Rayleigh.
Recent work of Choe \cite{C1,C2} shows that:\\
$\bullet$\ All sixth--root of unity
matroids are in fact ``strongly Rayleigh''  in a sense distinct from
the spectrum of conditions in Section 4.\\
Also:
$\bullet$\ A binary matroid is strongly Rayleigh if and only if it
is regular.\\
$\bullet$\ Every matroid with at most seven elements is Rayleigh.\\
$\bullet$\ Every matroid with a $2$--transitive automorphism group
is negatively correlated.\\
On the negative side:\\
$\bullet$\ There is a rank $4$ transversal matroid which is not
balanced.\\
In particular, such matroids need not be HPP, which settles negatively
a question left open in \cite{COSW}.\\
$\bullet$\ Every finite projective geometry fails to be HPP.\\
$\bullet$\ There is a balanced matroid which is not Rayleigh.\\
Combined with the results in Section 3, this shows that the class
of balanced matroids is not closed by taking $2$--sums.

We conclude in Section 6 with a few open problems.  For example:\\
$\bullet$\ Is every matroid of rank three a Rayleigh matroid?

We thank Jim Geelen, Criel Merino, Alan Sokal, and Dominic Welsh
for valuable converations on this subject, and Robert Shrock and
Earl Glen Whitehead, Jr. for invitation to a minisymposium on
``Graph Theory with Applications to Chemistry and Physics'' at the
First Joint Meeting of the C.A.I.M.S. and S.I.A.M. in Montreal, 
June 16--20, 2003, at which these results were presented.

\section{Balanced matroids.}

Feder and Mihail \cite{FM} prove two main results about balanced 
matroids.  First:\\
$\bullet$\ Every regular matroid is balanced.\\
This establishes
a large class of examples including, of course, all graphic or
cographic matroids.  (See Proposition 5.1 and Corollary 4.9 below.)
Second:\\
$\bullet$\ The basis--exchange graph of a balanced matroid
has cutset expansion at least one.\\
To explain this, the
\emph{basis--exchange graph} of a matroid $\M$ is the simple 
graph with the set of bases of $\M$ as its vertex--set, and with an
edge $B_{1}\sim B_{2}$ if and only if $|B_{1}\triangle B_{2}|=2$
(in which $\triangle$ denotes the symmetric difference of sets).  A simple
graph $G=(V,E)$ has \emph{cutset expansion at least} $\rho$ provided
that for every $\none\neq S\subset V$,
$$\frac{|\{e\in E:\ e\cap S\neq\varnothing\ \mathrm{and}\ 
e\cap(V\drop S)\neq\varnothing\}|}{\min\{|S|,|V\drop S|\}}\geq \rho.$$
Such isoperimetric inequalities imply that the natural
random walk on the graph converges rapidly to the uniform distribution
on the vertices.  This leads to an efficient algorithm
for generating a random basis of a balanced matroid approximately
uniformly. See \cite{FM} for details.

The matroid $\S_{8}$ is represented over $GF(2)$ by the matrix
$$\left[\begin{array}{llllllll}
1 & 1 & 1 & 1 & 1 & 1 & 1 & b\\
0 & 1 & 0 & 0 & 0 & 1 & 1 & 1\\
0 & 0 & 1 & 0 & 1 & 0 & 1 & 1\\
0 & 0 & 0 & 1 & 1 & 1 & 0 & 1
\end{array}\right]$$
with $b=0$, and the matroid $\A_{8}=\A\G(3,2)$ is represented over 
$GF(2)$ by this matrix with $b=1$.
Feder and Mihail refer to unpublished work showing that $\S_{8}$ is
the only minor--minimal binary non--balanced matroid.  To our 
knowledge, the only argument in print for this claim is in Chapter 5
and Appendix D of Merino's thesis \cite{Me}, but it contains an error.
Specifically, the argument rests on five points:\\
$\bullet$\ The matroid $\S_{8}$ is not negatively correlated.  This
was observed by Seymour and Welsh \cite{SW} and is not hard to 
verify.  (Labelling the ground--set $\{1,\ldots,8\}$ corresponding to
the columns of the above matrix, we have $(S_{8})_{1}=28$,
$(S_{8})_{8}=20$, $(S_{8})_{1,8}=12$, and $S_{8}=48$, so that
$\Delta S_{8}\{1,8\}=28\cdot 20-12\cdot 48=-16<0$.)\\
$\bullet$\ The matroid $\A_{8}$ is a ``splitter'' for the class
of binary matroids which do not contain an $\S_{8}$ minor.  More
explicitly, if a connected binary matroid $\M$ with no $\S_{8}$ minor
has $\A_{8}$ as a proper minor, then $\M$ can be expressed as a $2$--sum
with $\A_{8}$ as one of the factors.  This is an 
unpublished result of Seymour and is explained in Appendix D of 
\cite{Me}.\\
$\bullet$\ Every binary matroid which does not contain $\S_{8}$ or
$\A_{8}$ as a minor can be constructed from regular matroids,
the Fano matroid $\F_{7}$, and its dual $\F_{7}^{*}$ by taking direct
sums and $2$--sums.  This is due to Seymour \cite{Sey}.\\
$\bullet$\ The matroids $\A_{8}$, $\F_{7}$, and $\F_{7}^{*}$ are balanced.
This also is not  difficult to verify and appears in Appendix D of
\cite{Me}.\\
$\bullet$\ The class of balanced matroids is closed by taking 
$2$--sums.  This appears as Lemma 5.4.4 in \cite{Me}, but the
argument in support of it contains an error on the first part of page
$113$.  In fact, this claim is false (Theorem 5.11).

To explain the difficulty, consider a matroid $\M$ and distinct
elements $e,f,g$ of $E(\M)$.  Then, since $M=M_{g}+M^{g}$
\emph{et cetera}, a short calculation shows that
$$ \Delta M\{e,f\}=
\Delta M_{g}\{e,f\}+\Theta M\{e,f|g\}+\Delta M^{g}\{e,f\},$$
in which
\showon
\Delta M_{g}\{e,f\} &:=& M_{eg}^{f}M_{fg}^{e}-M_{efg}M_{g}^{ef},\\
\Delta M^{g}\{e,f\} &:=& M_{e}^{fg}M_{f}^{eg}-M_{ef}^{g}M^{efg},
\showoff
and
the \emph{central term for $\{e,f\}$ and $g$ in $\M$} is given by
$$\Theta M\{e,f|g\}:=M_{e}^{fg}M_{fg}^{e}+M_{f}^{eg}M_{eg}^{f}
-M_{g}^{ef}M_{ef}^{g}-M_{efg}M^{efg}.$$
Now let $\Q$ be another matroid, with $E(\Q)\cap E(\M)=\{g\}$, and
consider the $2$--sum  $\N=\M\oplus_{g}\Q$ of $\M$ and $\Q$ along
$g$.  The set of bases of $\N$ is
$$\N:=\{B_{1}\cup B_{2}:\ (B_{1},B_{2})\in
(\M_{g}\times \Q^{g})\cup(\M^{g}\times\Q_{g})\}$$
by definition, so that $N = M_{g}Q^{g}+M^{g}Q_{g}$.
Again, a short calculation shows that
\showon
& &\Delta N\{e,f\}=\\
& & (Q^{g})^{2}\Delta M_{g}\{e,f\}+Q^{g}Q_{g}\Theta M\{e,f|g\}+
(Q_{g})^{2}\Delta M^{g}\{e,f\}.
\showoff

Now assume that $\M$ is balanced.  If the class of
balanced matroids is closed by taking $2$--sums then
$\Delta N\{e,f\}\geq 0$ for any balanced choice of $\Q$.
That is, the quadratic polynomial
$$p(y):=y^{2}\Delta M_{g}\{e,f\}+y\Theta M\{e,f|g\}+
\Delta M^{g}\{e,f\}$$
is such that $p(\lambda)\geq 0$ for any real number of the form
$\lambda=Q^{g}/Q_{g}$ with $\Q$ balanced and $g\in E(\Q)$.

For positive integers $a$ and $b$, let $G(a,b)$ be the graph
obtained from a path with $b$ edges by replacing each edge by $a$
parallel edges, then joining the end--vertices by a new ``root'' edge.
Label the root edge of $G(a,b)$ by $g$.  The graphic (cycle) matroid
$\Q(a,b)$ of $G(a,b)$ is balanced by the result of Feder and Mihail.
Now, since $Q(a,b)^{g}/Q(a,b)_{g}=a/b$, every positive rational number
is of the form $\lambda$ above.

Therefore, the polynomial $p(y)$
above must satisfy $p(\lambda)\geq 0$ for all $\lambda\geq 0$, and
since both $\Delta M_{g}\{e,f\}$ and $\Delta M^{g}\{e,f\}$ are
nonnegative the zeros of $p(y)$ are either nonreal complex conjugates
or are real and of the same sign.  This implies that
$$\Theta M\{e,f|g\}\geq -2\sqrt{\Delta M_{g}\{e,f\}\Delta M^{g}\{e,f\}}.$$

This ``triple condition'' on the balanced matroid $\M$ is necessary 
for all $\{e,f\}$ and $g$ in $E(\M)$ if the class of balanced matroids
is closed by taking $2$--sums.  However, it is unclear 
whether or not this can be deduced from the hypothesis that
$\M$ is balanced.  The Rayleigh hypothesis, on the other hand, 
includes these triple conditions and can be carried through
the $2$--sum construction with ease, as we shall see in the next 
section.

\section{Rayleigh matroids.}

The term ``Rayleigh matroid'' is motivated by analogy with a property 
of electrical networks.  Consider a (multi)graph $G=(V,E)$ together with
a set $\y=\{y_{e}:\ e\in E\}$ of positive real numbers indexed by the
edges of $G$.  Thinking of each $y_{e}$ as the electrical conductance
of the edge $e\in E$, for any two vertices $a,b\in V$ we may ask for
the value of the effective conductance $\Y_{ab}(G;\y)$ of the graph
as a whole, considered as a network joining the poles $a$ and $b$.
In 1847, Kirchhoff \cite{Ki} proved that
$$\Y_{ab}(G;\y)=\frac{T(G;\y)}{T(G/ab;\y)},$$
in which $T(G;\y):=\sum_{T}\y^{T}$ with the sum over all spanning
trees of $G$, and $T(G/ab;\y)$ is defined similarly except that
$G/ab$ is the graph obtained from $G$ by merging $a$ and $b$ into a
single vertex.

It is physically intuitive that if $y_{c}>0$ for all $c\in E$ and
$y_{e}$ is increased, then $\Y_{ab}(G;\y)$ does not decrease -- this
property is called \emph{Rayleigh monotonicity}.  (This will be proven
below when we show that sixth--root of unity matroids -- in 
particular, graphs -- are Rayleigh matroids.)  Nonnegativity of 
$\partial\Y_{ab}(G;\y)/\partial y_{e}$ is equivalent to the 
inequality
$$\frac{\partial T(G;\y)}{\partial y_{e}}T(G/ab;\y)
\geq T(G;\y)\frac{\partial T(G/ab;\y)}{\partial y_{e}}.$$
Rephrasing this in terms of the graph $H$ obtained from $G$ by 
adjoining a new edge $f$ with ends $a$ and $b$, the inequality is
$$T_{e}^{f}(H;\y)T_{f}(H;\y) \geq T^{f}(H;\y)T_{ef}(H;\y),$$
in which $T_{e}^{f}(H;\y)$ is the sum of $\y^{T}$ over all spanning
trees $T$ of the graph obtained by contracting $e$ and deleting $f$
from $H$, \emph{et cetera}.  A little cancellation shows that this
is equivalent to the inequality
$$T_{e}^{f}(H;\y)T_{f}^{e}(H;\y)-T_{ef}(H;\y)T^{ef}(H;\y)\geq 0.$$
Replacing $T(H;\y)$ by the basis--generating polynomial $M(\y)$ of
a more general matroid $\M$, we arrive at the condition
$\Delta M\{e,f\}(\y)\geq 0$ defining Rayleigh matroids.

To simplify notation, when calculating with Rayleigh matroids we
will henceforth usually omit reference to the variables $\y$ -- writing
$M_{I}^{J}$ instead of $M_{I}^{J}(\y)$ \emph{et cetera} -- unless
a particular substitution of variables requires emphasis.  We will
also write ``$\y>\zero$'' as shorthand for ``$y_{c}>0$ for all $c\in 
E$'', and ``$\y\equiv\one$'' as shorthand for ``$y_{c}=1$ for all
$c\in E$''.

\begin{PROP}
A matroid $\M$ is Rayleigh if and only if the dual matroid $\M^{*}$
is Rayleigh.
\end{PROP}
\begin{proof}
For disjoint subsets $I,J\subseteq E$ we have $M_{I}^{*J}(\y)=\y^{E}
M_{J}^{I}(\one/\y)$, in which $\one/\y:=\{1/y_{c}:\ c\in E\}$.
Therefore, the inequality $\Delta M^{*}\{e,f\}(\y)\geq 0$
is equivalent to the inequality $\Delta M\{e,f\}(\one/\y)\geq 0$.
From this the result follows.
\end{proof}

\begin{PROP}
If $\M$ is a Rayleigh matroid and $\N$ is a minor of $\M$ then $\N$ is
a Rayleigh matroid.
\end{PROP}
\begin{proof}
Since $\M$ is Rayleigh, for distinct $e,f,g\in E$ and $\y>\zero$
we have
$$\Delta M\{e,f\}=
y_{g}^{2}\Delta M_{g}\{e,f\}+y_{g}\Theta M\{e,f|g\}+\Delta 
M^{g}\{e,f\}\geq 0.$$
Take the limit of this as $y_{g}\goesto 0$ to see that
$\Delta M^{g}\{e,f\}\geq 0$.  Since $e,f\in E(\M^{g})$
and $\y>\zero$ are arbitrary, this shows that $\M^{g}$ is Rayleigh.
Similarly, by considering the limit of $y_{g}^{-2}\Delta M\{e,f\}$
as $y_{g}\goesto \infty$ we see that $\M_{g}$ is Rayleigh.
The case of a general minor is obtained by iteration of the 
above two cases.
\end{proof}

\begin{CORO}
Every Rayleigh matroid $\M$ is balanced and satisfies the triple
condition
$$\Theta M\{e,f|g\}\geq -2\sqrt{\Delta M_{g}\{e,f\}\Delta M^{g}\{e,f\}}$$
for distinct $e,f,g\in E(\M)$ when $\y>\zero$.
\end{CORO}
\begin{proof}
If $\M$ is a Rayleigh matroid then by setting $\y\equiv\one$
we see that $\M$ is negatively correlated.  Since every
minor of $\M$ is also Rayleigh, it follows that $\M$ is balanced.

For distinct $e,f,g\in E(\M)$, when $y_{c}>0$ for all $c\neq g$,
the polynomial
$$\Delta M\{e,f\}=
y_{g}^{2}\Delta M_{g}\{e,f\}+y_{g}\Theta M\{e,f|g\}+\Delta 
M^{g}\{e,f\}$$
in $y_{g}$ is nonnegative for all $y_{g}>0$.   As in Section 2,
this implies the desired inequality.
\end{proof}

\begin{PROP}
Let $\M$ be a matroid with ground set $E$, and let $I,J$ be disjoint
subsets of $E$.  If $\M$ is Rayleigh and $\y>\zero$
then $$M_{I}M_{J}\geq M_{IJ}M.$$
\end{PROP}
\begin{proof}
The inequality is trivial if either $I$ or $J$ is dependent, so assume
that both $I$ and $J$ are independent in $\M$.

We first prove the result for $I=\{e_{1}\}$ and
$J=\{f_{1},\ldots,f_{k}\}$.
Notice that the Rayleigh difference of $\{e,f\}$ in $\M$ may
also be expressed as $\Delta M\{e,f\}=M_{e}M_{f}-M_{ef}M$.  Thus,
the Rayleigh condition is that if $\y>\zero$ then
$M_{e}M_{f}\geq M_{ef}M$.  Since every (contraction) minor of $\M$
is also Rayleigh, we see that if $\y>\zero$ then
$$\frac{M_{e_{1}}}{M}\geq\frac{M_{e_{1}f_{1}}}{M_{f_{1}}}\geq
\frac{M_{e_{1}f_{1}f_{2}}}{M_{f_{1}f_{2}}}\geq\cdots\geq
\frac{M_{e_{1}J}}{M_{J}}.$$
That is, $M_{I}M_{J}\geq M_{IJ}M$ in this case.

Viewed another way, we have shown that if $\M$ is Rayleigh and
$\y>\zero$ then for any non--loop $e_{1}\in E$ and $J\subseteq E$,
$M_{J}/M\geq M_{e_{1}J}/M_{e_{1}}$.  If now $I=\{e_{1},e_{2},\ldots,
e_{m}\}$ is independent then since each (contraction) minor of $\M$
is Rayleigh
$$\frac{M_{J}}{M}\geq\frac{M_{e_{1}J}}{M_{e_{1}}}
\geq\frac{M_{e_{1}e_{2}J}}{M_{e_{1}e_{2}}}
\geq\cdots\geq\frac{M_{IJ}}{M_{I}}.$$
This implies the desired inequality.
\end{proof}
The probability space associated with $\M$ and $\y>\zero$
assigns to each basis $B$ of $\M$ the probability $\y^{B}/M(\y)$.
As in \cite{FM,Ly}, Proposition 3.4 leads to the fact that any two
increasing events with disjoint support in this space are negatively 
correlated, provided that $\M$ is Rayleigh.

\begin{THM}
Let $\M$ and $\Q$ be matroids with $E(\M)\cap E(\Q)=\{g\}$, and let
$\N = \M\oplus_{g}\Q$ be the $2$--sum of $\M$ and $\Q$ along $g$.
If $\M$ and $\Q$ are Rayleigh matroids then $\N$ is a Rayleigh matroid.
\end{THM}
\begin{proof}
Fix $y_{c}>0$ for all $c\in E(\N)$, and consider any $e,f\in E(\N)$.
We must show that $\Delta N\{e,f\}\geq 0$.  Up to symmetry of the 
hypotheses there are essentially two cases:\\
(i)\ $e\in E(\M)\drop\{g\}$ and $f\in E(\Q)\drop\{g\}$;\\
(ii)\ $\{e,f\}\subseteq E(\M)\drop\{g\}$.

For case (i) a short calculation using $N=M_{g}Q^{g}+M^{g}Q_{g}$
\emph{et cetera} shows that
$$\Delta N\{e,f\}=\Delta M\{e,g\}\Delta Q\{f,g\}.$$
Since $\M$ and $\Q$ are Rayleigh and $\y>\zero$, both factors on the right are
nonnegative, so that $\Delta N\{e,f\}\geq 0$ as well.

For case (ii) we calculate that
\showon
& &\Delta N\{e,f\}(\y)=\\
& & (Q^{g})^{2}\Delta M_{g}\{e,f\}+Q^{g}Q_{g}\Theta M\{e,f|g\}+
(Q_{g})^{2}\Delta M^{g}\{e,f\}.
\showoff
If $Q^{g}(\y)=0$ or $Q_{g}(\y)=0$ then $\Delta N\{e,f\}\geq 0$
because both $\M^{g}$ and $\M_{g}$ are Rayleigh.
Otherwise, by defining $w_{c}:=y_{c}$ for all
$c\in E(\M)\drop\{g\}$ and $w_{g}:=Q^{g}(\y)/Q_{g}(\y)$, we see
that
$$\Delta N\{e,f\}(\y)=(Q_{g})^{2}\Delta M\{e,f\}(\w)\geq 0,$$
since $\w>\zero$ and $\M$ is Rayleigh.

This proves that $\N=\M\oplus_{g}\Q$ is Rayleigh.
\end{proof}

For a matroid $\M$ and a set $\m:=\{m_{e}:\ e\in E(\M)\}$ of
positive integers indexed by $E(\M)$, let $\M[\m]$ be the
matroid obtained from $\M$ by replacing each element $e\in E(\M)$
by a parallel class of $m_{e}$ elements.  Equivalently, $\M[\m]$
is obtained from $\M$ by attaching the uniform matroid 
$\U_{1,1+m_{e}}$ to $\M$ by a $2$--sum along $e$, for each $e\in 
E(\M)$.

\begin{THM}
For a matroid $\M$, the following conditions are equivalent:\\
\textup{(a)}\ The matroid $\M$ is Rayleigh;\\
\textup{(b)}\ Every matroid of the form $\M[\m]$ is Rayleigh;\\
\textup{(c)}\ Every matroid of the form $\M[\m]$ is balanced;\\
\textup{(d)}\ Every matroid of the form $\M[\m]$ is negatively
correlated.
\end{THM}
\begin{proof}
To see that (a) implies (b) we note that a matroid $\Q$ of rank one
is Rayleigh since for any $e,f\in E(\Q)$ we have $\Q_{ef}\equiv 0$.
Thus, since $\M[\m]$ is expressed as a $2$--sum of Rayleigh matroids
it is also Rayleigh, by Theorem 3.5.

That (b) implies (c) is immediate from Corollary 3.3, and (c)
implies (d) is immediate from the definitions.

To see that (d) implies (a) assume that $\M$ is not Rayleigh.
Thus, there exist distinct $e,f\in E(\M)$ and positive real numbers
$\y>\zero$ such that $\Delta M\{e,f\}(\y) <0$.   Since the rational
numbers are dense in the real numbers, there are positive rationals
$\q=\{q_{c}:\ c\in E\}$ such that $\Delta M\{e,f\}(\q) <0$.
Let $D$ be the smallest positive common denominator of all the numbers
$\{q_{c}:\ c\in \drop\{e,f\}\}$, and for $c\in E\drop\{e,f\}$ let
$m_{c}:=Dq_{c}$, a positive integer.   Since $\Delta\M\{e,f\}(\y)$
is independent of $y_{e}$ and $y_{f}$ we may put $m_{e}:=m_{f}:=1$.
Since  $\Delta\M\{e,f\}(\y)$ is homogeneous of degree $2r-2$ (where $r$
is the rank of $\M$) we have
$$\Delta\M\{e,f\}(\m)=D^{2r-2}\Delta\M\{e,f\}(\q)<0.$$
However, we also have $\Delta\M[\m]\{e,f\}(\one)=\Delta\M\{e,f\}(\m)<0$,
so that $\M[\m]$ is not negatively correlated.
\end{proof}

\begin{CORO}
The following statements are equivalent:\\
\textup{(a)}\ Every balanced matroid is Rayleigh;\\
\textup{(b)}\ The class of balanced matroids is closed by taking
$2$--sums.
\end{CORO}
\begin{proof}
To show that (a) implies (b), let $\M$ and $\Q$ be balanced matroids
such that $E(\M)\cap E(\Q)=\{g\}$.  By (a) both $\M$ and $\Q$ are
Rayleigh, so that $\M\oplus_{g}\Q$ is
Rayleigh by Theorem 3.5, and hence balanced by Corollary 3.3.

To show that (b) implies (a), consider a balanced matroid $\M$.
Since uniform matroids of rank one are balanced, the hypothesis (b)
implies that every matroid of the form $\M[\m]$ is balanced.
By Theorem 3.6, it follows that $\M$ is Rayleigh.
\end{proof}
\noindent
In Theorem 5.11 we will see that the two statements of Corollary 3.7
are in fact false.

\begin{THM}
A binary matroid is Rayleigh if and only if it does not contain 
$\S_{8}$ as a minor.
\end{THM}
\begin{proof}
The outline of the argument has been sketched in Section 2 (for 
balanced matroids in place of Rayleigh matroids).
For the first point, since $\S_{8}$ is not negatively correlated it
is not balanced, hence not Rayleigh.  The second and third points need no
revision, and the fifth point is substantiated for Rayleigh matroids
by Theorem 3.5.

It remains to show that the matroids $\A_{8}$, 
$\F_{7}$, and $\F_{7}^{*}$ are Rayleigh.  Since $\F_{7}$ is obtained 
from $\A_{8}$ by contracting any element, Propositions 3.1 and 3.2 
imply that it is enough to show that $\A_{8}$ is Rayleigh.
Let the ground--set of $\A_{8}$ be $E=\{1,\ldots,8\}$
corresponding to the columns of the representing matrix in Section 2.
The automorphism group of $\A_{8}$ is $2$--transitive on $E$, so
in order to check that this matroid is Rayleigh it suffices to show that
$\Delta A_{8}\{7,8\}\geq 0$ when $\y>\zero$.
A direct computation with the aid of \textsc{Maple 6.01} shows that
\showon
& & \Delta A_{8}\{7,8\}\\
&=&
2\,{y_{{1}}}^{2}{y_{{2}}}^{2}y_{{5}}y_{{6}}+
2\,{y_{{1}}}^{2}y_{{2}}y_{{3}}y_{{4}}y_{{6}}+
2\,{y_{{1}}}^{2}y_{{2}}y_{{3}}y_{{5}}y_{{6}}+
2\,{y_{{1}}}^{2}y_{{2}}y_{{3}}{y_{{6}}}^{2}\\
&+&
2\,{y_{{1}}}^{2}y_{{2}}y_{{4}}y_{{5}}y_{{6}}+
2\,{y_{{1}}}^{2}y_{{2}}y_{{4}}{y_{{6}}}^{2}+
2\,{y_{{1}}}^{2}y_{{2}}{y_{{5}}}^{2}y_{{6}}+
2\,{y_{{1}}}^{2}y_{{2}}y_{{5}}{y_{{6}}}^{2}\\
&+&
2\,{y_{{1}}}^{2}{y_{{3}}}^{2}y_{{4}}y_{{6}}+
2\,{y_{{1}}}^{2}y_{{3}}{y_{{4}}}^{2}y_{{6}}+
2\,{y_{{1}}}^{2}y_{{3}}y_{{4}}y_{{5}}y_{{6}}+
2\,{y_{{1}}}^{2}y_{{3}}y_{{4}}{y_{{6}}}^{2}\\
&+&
2\,{y_{{1}}}^{2}y_{{3}}y_{{5}}{y_{{6}}}^{2}+
2\,{y_{{1}}}^{2}y_{{4}}y_{{5}}{y_{{6}}}^{2}+
2\,y_{{1}}{y_{{2}}}^{2}y_{{3}}y_{{4}}y_{{5}}+
2\,y_{{1}}{y_{{2}}}^{2}y_{{3}}{y_{{5}}}^{2}\\
&+&
2\,y_{{1}}{y_{{2}}}^{2}y_{{3}}y_{{5}}y_{{6}}+
2\,y_{{1}}{y_{{2}}}^{2}y_{{4}}{y_{{5}}}^{2}+
2\,y_{{1}}{y_{{2}}}^{2}y_{{4}}y_{{5}}y_{{6}}+
2\,y_{{1}}{y_{{2}}}^{2}{y_{{5}}}^{2}y_{{6}}\\
&+&
2\,y_{{1}}{y_{{2}}}^{2}y_{{5}}{y_{{6}}}^{2}+
2\,y_{{1}}y_{{2}}{y_{{3}}}^{2}{y_{{4}}}^{2}+
2\,y_{{1}}y_{{2}}{y_{{3}}}^{2}y_{{4}}y_{{5}}+
2\,y_{{1}}y_{{2}}{y_{{3}}}^{2}y_{{4}}y_{{6}}\\
&+&
2\,y_{{1}}y_{{2}}y_{{3}}{y_{{4}}}^{2}y_{{5}}+
2\,y_{{1}}y_{{2}}y_{{3}}{y_{{4}}}^{2}y_{{6}}+
2\,y_{{1}}y_{{2}}y_{{3}}y_{{4}}{y_{{5}}}^{2}+
4\,y_{{1}}y_{{2}}y_{{3}}y_{{4}}y_{{5}}y_{{6}}\\
&+&
2\,y_{{1}}y_{{2}}y_{{3}}y_{{4}}{y_{{6}}}^{2}+
2\,y_{{1}}y_{{2}}y_{{3}}{y_{{5}}}^{2}y_{{6}}+
2\,y_{{1}}y_{{2}}y_{{3}}y_{{5}}{y_{{6}}}^{2}+
2\,y_{{1}}y_{{2}}y_{{4}}{y_{{5}}}^{2}y_{{6}}\\
&+&
2\,y_{{1}}y_{{2}}y_{{4}}y_{{5}}{y_{{6}}}^{2}+
2\,y_{{1}}y_{{2}}{y_{{5}}}^{2}{y_{{6}}}^{2}+
2\,y_{{1}}{y_{{3}}}^{2}{y_{{4}}}^{2}y_{{5}}+
2\,y_{{1}}{y_{{3}}}^{2}{y_{{4}}}^{2}y_{{6}}\\
&+&
2\,y_{{1}}{y_{{3}}}^{2}y_{{4}}y_{{5}}y_{{6}}+
2\,y_{{1}}{y_{{3}}}^{2}y_{{4}}{y_{{6}}}^{2}+
2\,y_{{1}}y_{{3}}{y_{{4}}}^{2}y_{{5}}y_{{6}}+
2\,y_{{1}}y_{{3}}{y_{{4}}}^{2}{y_{{6}}}^{2}\\
&+&
2\,y_{{1}}y_{{3}}y_{{4}}y_{{5}}{y_{{6}}}^{2}+
2\,{y_{{2}}}^{2}{y_{{3}}}^{2}y_{{4}}y_{{5}}+
2\,{y_{{2}}}^{2}y_{{3}}{y_{{4}}}^{2}y_{{5}}+
2\,{y_{{2}}}^{2}y_{{3}}y_{{4}}{y_{{5}}}^{2}\\
&+&
2\,{y_{{2}}}^{2}y_{{3}}y_{{4}}y_{{5}}y_{{6}}+
2\,{y_{{2}}}^{2}y_{{3}}{y_{{5}}}^{2}y_{{6}}+
2\,{y_{{2}}}^{2}y_{{4}}{y_{{5}}}^{2}y_{{6}}+
2\,y_{{2}}{y_{{3}}}^{2}{y_{{4}}}^{2}y_{{5}}\\
&+&
2\,y_{{2}}{y_{{3}}}^{2}{y_{{4}}}^{2}y_{{6}}+
2\,y_{{2}}{y_{{3}}}^{2}y_{{4}}{y_{{5}}}^{2}+
2\,y_{{2}}{y_{{3}}}^{2}y_{{4}}y_{{5}}y_{{6}}+
2\,y_{{2}}y_{{3}}{y_{{4}}}^{2}{y_{{5}}}^{2}\\
&+&
2\,y_{{2}}y_{{3}}{y_{{4}}}^{2}y_{{5}}y_{{6}}+
2\,y_{{2}}y_{{3}}y_{{4}}{y_{{5}}}^{2}y_{{6}}+
2\,{y_{{3}}}^{2}{y_{{4}}}^{2}y_{{5}}y_{{6}}\\
&+&
\left (y_{{1}}y_{{5}}y_{{6}}-y_{{3}}y_{{4}}y_{{5}}\right )^{2}+
\left (y_{{1}}y_{{2}}y_{{5}}-y_{{1}}y_{{3}}y_{{4}}\right )^{2}+
\left (y_{{2}}y_{{4}}y_{{5}}-y_{{1}}y_{{4}}y_{{6}}\right )^{2}\\
&+&
\left (y_{{2}}y_{{3}}y_{{5}}-y_{{1}}y_{{3}}y_{{6}}\right )^{2}+
\left (y_{{2}}y_{{5}}y_{{6}}-y_{{3}}y_{{4}}y_{{6}}\right )^{2}+
\left (y_{{1}}y_{{2}}y_{{6}}-y_{{2}}y_{{3}}y_{{4}}\right )^{2}
\showoff
Since this is clearly nonnegative for $\y>\zero$ we see that
$\A_{8}$ is Rayleigh.  This completes the proof.
\end{proof}

\begin{CORO}
A binary matroid is balanced if and only if it is Rayleigh.
\end{CORO}
\begin{proof}
By Corollary 3.3, every Rayleigh matroid is balanced.
If $\M$ is a balanced matroid then $\M$ does not contain $\S_{8}$
as a minor, since $\S_{8}$ is not negatively correlated.  If
$\M$ is also binary then $\M$ is Rayleigh, by Theorem 3.8.
\end{proof}

\section{Half--plane property matroids.}

A polynomial $P(\y)=\sum_{\alpha}c_{\alpha}\y^{\alpha}$ in several
complex variables $\y=\{y_{e}:\ e\in E\}$ has the \emph{half--plane
property} provided that whenever Re$(y_{e})>0$ for all $e\in E$,
then $P(\y)\neq 0$.  We say that a matroid $\M=(E,\B)$ is a
\emph{half--plane property matroid} (HPP matroid, for short) if its
basis--generating
polynomial $M(\y):=\sum_{B\in\B}\y^{B}$ has the half--plane property.
This class of 
polynomials is investigated thoroughly in \cite{COSW}, from which we
take the following facts without proof.

\begin{LMA}[\cite{COSW}, Proposition $4.2$]
Let $P(\y)$ be a polynomial in the variables
$\y=\{y_{e}:\ e\in E\}$, and let $d_{e}$ be the degree of $y_{e}$ in
$P$ for each $e\in E$.  If $P(\y)$ has the half--plane property then
$\y^{\d}P(\one/\y)$ has the half--plane property.
\end{LMA}

\begin{LMA}[\cite{FB}, Theorem $18$, or \cite{COSW}, Proposition $3.4$.]
Let $P(\y)$ be a polynomial in the variables
$\y=\{y_{e}:\ e\in E\}$, fix $e\in E$, and let
$P(\y)=\sum_{j=0}^{n}P_{j}(y_{c}:\ c\neq e\})y_{e}^{j}$.
If $P$ has the half--plane property then each $P_{j}$ has the 
half--plane property.
\end{LMA}

\begin{LMA}[\cite{COSW}, Proposition $5.2$]
Let $P(\y)$ be a homogeneous polynomial in the variables
$\y=\{y_{e}:\ e\in E\}$.  For nonnegative real numbers
$\a=\{a_{e}:\ e\in E\}$ and $\b=\{b_{e}:\ e\in E\}$, let
$P(\a x+\b)$ be the polyomial obtained by substituting
$y_{e}=a_{e}x+b_{e}$ for all $e\in E$.  The following are
equivalent:\\
\textup{(a)} $P(\y)$ has the half--plane property;\\
\textup{(b)} for all sets of nonnegative real numbers $\a$ and $\b$,
$P(\a x + \b)$ has only real zeros.
\end{LMA}

\begin{PROP}[\cite{COSW}, Propositions $3.1$, $4.1$, and $4.2$]
The class of HPP matroids is closed by taking duals and minors.
\end{PROP}
\begin{proof}
For a matroid $\M$ on a set $E$, the dual matroid $\M^{*}$ has
basis generating polynomial $M^{*}(\y)=\y^{E}M(\one/\y)$.  By Lemma
4.1, if $\M$ is HPP then $\M^{*}$ is HPP.  For $g\in E$ we have
$M(\y)=y_{g}M_{g}(\y)+M^{g}(\y)$.  Lemma 4.2 implies that if $\M$ is
HPP then both $\M_{g}$ and $\M^{g}$ are HPP.  The case of a general
minor of $\M$ follows by iterating these two cases.
\end{proof}
Many other operations are shown to preserve the half--plane property
in Section $4$ of \cite{COSW}, $2$--sums in particular.

Theorem 4.5 was proven for regular matroids and $\y\equiv\one$
by Godsil \cite{Go}.
\begin{THM}
Let $\M$ be a matroid on a set $E$.  Let 
$(S,T,C_{1},\ldots,C_{k})$ be an ordered partition of $E$ into
pairwise disjoint nonempty subsets, and fix nonnegative integers
$c_{1},\ldots,c_{k}$.  For each $0\leq j\leq |S|$, let
$M_{j}(\y):=\sum_{B}\y^{B}$, with the sum over all bases $B$
of $\M$ such that $|B\cap S|=j$ and $|B\cap C_{i}|=c_{i}$ for all
$1\leq i\leq k$.  If $\M$ is a HPP matroid and $\y>\zero$,
then the polynomial $\sum_{j=0}^{|S|}M_{j}(\y)x^{j}$ in the
variable $x$ has only real zeros.
\end{THM}
\begin{proof}
Let $\M$ be a HPP matroid and fix $\y>\zero$.
Let $s$, $t$, and $z_{1},\ldots,z_{k}$ be 
indeterminates, and for $e\in E$ put
$$u_{e}:=\left\{\begin{array}{ll}
y_{e}s & \mathrm{if}\ \  e\in S,\\
y_{e}t & \mathrm{if}\ \  e\in T,\\
y_{e}z_{i} & \mathrm{if}\ \ e\in C_{i}.\end{array}\right.$$
Then $M(\u)$ is a homogeneous polynomial with the half--plane property
in the variables $s,t,z_{1},\ldots,z_{k}$.  By repeated application
of Lemma 4.2, the coefficient $M_{\c}(s,t)$ of
$z_{1}^{c_{1}}\cdots z_{k}^{c_{k}}$ in
$M(\u)$ also has the half--plane property, and is homogeneous.  In fact,
$$M_{\c}(s,t) = \sum_{j=0}^{|S|}M_{j}(\y)s^{j}t^{d-j},$$
in which $d=\mathrm{rank}(\M)-(c_{1}+\cdots+c_{k})$.
Upon substituting $s=x$ and $t=1$ in $M_{\c}(s,t)$, Lemma 4.3
implies that $\sum_{j=0}^{|S|}M_{j}(\y)x^{j}$ has only real zeros,
as claimed.
\end{proof}

Newton's Inequalities (item (51) of \cite{HLP}) state that
if a polynomial $\sum_{j=0}^{n}a_{j}x^{j}$ with real coefficients
has only real zeros then $\binom{n}{j}^{-2}a_{j}^{2}\geq\binom{n}{j-1}^{-1}
\binom{n}{j+1}^{-1}a_{j-1}a_{j+1}$ for all $1\leq j\leq n-1$.
That is, the sequence $\{\binom{n}{j}^{-1}a_{j}\}$ is 
\emph{logarithmically concave}.  Thus, Theorem 4.5 implies the
following corollary, first proved for regular matroids and 
$\y\equiv\one$ by Stanley \cite{St}.
\begin{CORO}
With the hypothesis and notation of Theorem $4.5$, for
each $1\leq j\leq |S|-1$,
$$\frac{M_{j}(\y)^{2}}{\binom{|S|}{j}^{2}}\geq
\frac{M_{j-1}(\y)}{\binom{|S|}{j-1}}\cdot
\frac{M_{j+1}(\y)}{\binom{|S|}{j+1}}.$$
\end{CORO}
\noindent
Corollary 4.6 can be viewed as a quantitative strengthening of the
basis exchange axiom for HPP matroids, as requested in Question 13.9
of \cite{COSW}.

For a subset $S\subseteq E(\M)$ and natural number $j$, let
$M(S,j;\y)=\sum_{B}\y^{B}$, with the sum over all bases $B$ of $\M$
such that $|B\cap S|=j$.  For each positive integer $m$, consider the
following conditions on a matroid $\M$:\\

\noindent
\textbf{RZ[$m$]:}\ \ If $\y>\zero$ then for all
$S\subseteq E$ with $|S|\leq m$ the polynomial $\sum_{j=0}^{|S|}
M(S,j;\y)x^{j}$ has only real zeros.\\

\noindent
\textbf{LC[$m$]:}\ \  If $\y>\zero$ then for all
$S\subseteq E$ with $|S|\leq m$ the sequence $\{\binom{|S|}{j}^{-1}
M(S,j;\y)\}$ is logarithmically concave.\\

\noindent
The $k=0$ case of Theorem 4.5 implies that a HPP matroid is
RZ[$m$] for all $m$, and Newton's Inequalities show that RZ[$m$]
implies LC[$m$] for every $m$.  The implications RZ[$m$]\
$\Longrightarrow$ RZ[$m-1$] and LC[$m$]\ $\Longrightarrow$ LC[$m-1$]
are trivial, as are the conditions RZ[$1$] and LC[$1$].  Thus, the
weakest nontrivial condition among these is LC[$2$].

\begin{LMA}
Let $\M$ be a matroid on the set $E$.  If $\M$ is Rayleigh and
$\y>\zero$ then for any $S\subseteq E$ with $|S|\geq 2$,
$$\frac{M(S,1;\y)^{2}}{|S|^{2}}\geq \binom{|S|}{2}^{-1}M(S,0;\y)M(S,2;\y).$$
\end{LMA}
\begin{proof}
For any real numbers $R_{1},\ldots,R_{m}$ with $m\geq 2$,
\showon
(R_{1}+\cdots+R_{m})^{2} &=& \sum_{i=1}^{m}\sum_{j=1}^{m}R_{i}R_{j}\\
&=&\sum_{\{i,j\}\subseteq\{1,\ldots,m\}}\left(2R_{i}R_{j}+
\frac{R_{i}^{2}+R_{j}^{2}}{m-1}\right)\\
&\geq& \frac{2m}{m-1}\sum_{\{i,j\}\subseteq\{1,\ldots,m\}}R_{i}R_{j},
\showoff
since $R_{i}^{2}+R_{j}^{2}\geq 2R_{i}R_{j}$.  Apply this inequality
when $S=\{e_{1},\ldots,e_{m}\}$ and $R_{i}:=y_{e_{i}}M_{e_{i}}^{S\drop 
e_{i}}(\y)$ for $1\leq i\leq m$, with the result that
\showon
M(S,1;\y)^{2}
&\geq& \frac{2|S|}{|S|-1}\sum_{\{e,f\}\subseteq S}
y_{e}y_{f}M_{e}^{S\drop e}(\y)M_{f}^{S\drop f}(\y)\\
&\geq& \frac{2|S|}{|S|-1}\sum_{\{e,f\}\subseteq S}
y_{e}y_{f}M_{ef}^{S\drop ef}(\y)M^{S}(\y)\\
&=& \frac{2|S|}{|S|-1}M(S,0;\y)M(S,2;\y).
\showoff
The second inequality uses the fact that each of the deletion minors 
$\M^{S\drop ef}$ of $\M$ is Rayleigh.  This is equivalent to the
stated inequality.
\end{proof}

\begin{THM}
The following conditions are equivalent:\\
\textup{(a)}\ \ the matroid $\M$ is \textup{LC[$2$]};\\
\textup{(b)}\ \ the matroid $\M$ is \textup{RZ[$2$]};\\
\textup{(c)}\ \ the matroid $\M$ is Rayleigh;\\
\textup{(d)}\ \ the matroid $\M$ is \textup{LC[$3$]}.
\end{THM}
\begin{proof}
Conditions (a) and (b) are equivalent because a quadratic polynomial
has only real zeros if and only if its discriminant is nonnegative.

To show that (a) implies (c) assume that $\M$ is LC[$2$],
and choose distinct $e,f\in E$.  Since $\M$ is LC[$2$],
if $w_{c}>0$ for all $c\in E$ then
$$\left(w_{e}M_{e}^{f}(\w)+w_{f}M_{f}^{e}(\w)\right)^{2}\geq
4w_{e}w_{f}M_{ef}(\w)M^{ef}(\w).$$
In particular, if $\y>\zero$ then let
$$w_{c}:=\left\{\begin{array}{ll}
y_{c} & \mathrm{if}\ \ c\not\in\{e,f\},\\
M_{f}^{e}(\y) & \mathrm{if}\ \ c=e,\\
M_{e}^{f}(\y) & \mathrm{if}\ \ c=f.\end{array}\right.$$
The inequality above becomes
$$\left(2M_{e}^{f}(\y)M_{f}^{e}(\y)\right)^{2}
\geq 4M_{e}^{f}(\y)M_{f}^{e}(\y)M_{ef}(\y)M^{ef}(\y).$$
After some cancellation, this shows that
$$M_{e}^{f}(\y)M_{f}^{e}(\y)\geq M_{ef}(\y)M^{ef}(\y).$$
Hence, $\M$ is Rayleigh.

To show that (c) implies (d) assume that $\M$ is Rayleigh, and
let $\y>\zero$.  For a subset $S\subseteq E$ with $|S|\geq 2$,
Lemma 4.7 shows that $|S|^{-2}M(S,1;\y)^{2}\geq\binom{|S|}{2}^{-1}
M(S,0;\y)M(S,2;\y)$.  This implies that $\M$ is LC[$2$] and verifies
one of the inequalities of the condition LC[$3$] when $|S|=3$. It
remains to show that if $|S|=3$ then  $M(S,2;\y)^{2}\geq 
3M(S,1;\y)M(S,3;\y)$.  To do this we apply Lemma 4.7 to the dual
matroid $\M^{*}$, which is also Rayleigh.  Since
$$M^{*}(S,j;\y)=\y^{E}M(S,3-j;\one/\y)$$ for $0\leq j\leq 3$,
Lemma 4.7 implies the required inequality,
showing that $\M$ is LC[$3$].

That (d) implies (a) is trivial.  This completes the proof.
\end{proof}

\begin{CORO}
Every HPP matroid is a Rayleigh matroid.
\end{CORO}
\begin{proof}
By Theorem 4.5, every HPP matroid satisfies the conditions RZ[$m$]
for all $m$; in particular, it satisfies RZ[$2$] and hence is
Rayleigh by Theorem 4.8.
\end{proof}

\section{Examples.}

A matrix $A$ of complex numbers is a \emph{sixth--root of unity}
matrix provided that every nonzero minor of $A$ is a sixth--root
of unity.  A matroid $\M$ is a \emph{sixth--root of unity matroid}
provided that it can be represented over the complex numbers by a
sixth--root of unity matrix.  For example, every regular matroid
is a sixth--root of unity matroid.  Whittle \cite{Wh} has shown that
a matroid is a sixth--root of unity matroid if and only if it is
representable over both $GF(3)$ and $GF(4)$.  For graphs, Proposition $5.1$
is part of the ``folklore'' of electrical engineering.  We take it
from Corollary $8.2(a)$ and Theorem $8.9$ of \cite{COSW}, but include
the short and interesting proof for completeness.

\begin{PROP}
Every sixth--root of unity matroid is a HPP matroid.
\end{PROP}
\begin{proof}
Let $A$ be a sixth--root of unity matrix of full row--rank $r$, 
representing the matroid $\M$, and let $A^{*}$ denote the 
conjugate transpose of $A$.  Index the columns of $A$ by
the set $E$, and let $Y:=\mathrm{diag}(y_{e}:\ e\in E)$ be a
diagonal matrix of indeterminates.  For an $r$--element subset $S
\subseteq E$, let $A[S]$ denote the square submatrix of $A$ supported on
the set $S$ of columns.  By the Binet--Cauchy formula,
$$\det(AYA^{*})=\sum_{S\subseteq E:\ |S|=r}|\det A[S]|^{2}\y^{S}=M(\y)$$
is the basis--generating polynomial of $\M$,  since $|\det A[S]|^{2}$
is $1$ or $0$ according to whether or not $S$ is a basis of $\M$.

Now we claim that if Re$(y_{e})>0$ for all $e\in E$, then $AYA^{*}$ is
nonsingular.  This suffices to prove the result.  Consider any nonzero
vector $\v\in\CC^{r}$.  Then $A^{*}\v\neq\zero$ since the columns of
$A^{*}$ are linearly independent.  Therefore
$$\v^{*}AYA^{*}\v=\sum_{e\in E}y_{e}|(A^{*}\v)_{e}|^{2}$$
has strictly positive real part, since for all $e\in E$ the numbers
$|(A^{*}\v)_{e}|^{2}$ are nonnegative reals and at least one of these
is positive.  In particular, for any nonzero $\v\in\CC^{r}$, the
vector $AYA^{*}\v$ is nonzero.  It follows that $AYA^{*}$ is 
nonsingular, completing the proof.
\end{proof}
The same proof shows that for any complex matrix $A$ of full row--rank
$r$, the polynomial
$$\det(AYA^{*})=\sum_{S\subseteq E:\ |S|=r}|\det A[S]|^{2}\y^{S}$$
has the half--plane property.  The weighted analogue of Rayleigh 
monotonicity in this case is discussed from a probabilistic point of 
view by Lyons \cite{Ly}.  It is a surprising fact that a complex
matrix $A$ of full row--rank $r$ has $|\det A[S]|^{2}=1$ for all
nonzero rank $r$ minors if and only if $A$ represents a sixth--root of
unity matroid (Theorem $8.9$ of \cite{COSW}).

Regarding converses to Proposition 5.1, we note the 
following:\\
$\bullet$\ A binary matroid is HPP if and
only if it is regular (Corollary 8.16 of \cite{COSW}).\\
$\bullet$\ A ternary matroid is HPP if and only if it
is a sixth--root of unity matroid (Corollary 8.17 of \cite{COSW}).\\
$\bullet$\ Every matroid representable over $GF(4)$ which is
shown to be HPP in \cite{COSW} is a sixth--root of unity matroid.
However, some unsettled cases are expected to be HPP
but not sixth--root of unity.\\
$\bullet$\ Every uniform matroid is HPP (Theorem 9.1 of \cite{COSW}).

Another class of examples of HPP matroids can be produced
using the Heilmann--Lieb Theorem (Theorem $4.6$ and Lemma $4.7$ of
\cite{HL}, or Theorem $10.1$ of \cite{COSW}), but we have nothing new
to add here.

Proposition 5.1 and Corollary 4.9 show that every sixth--root of
unity matroid is Rayleigh.  This implies the result of Feder and
Mihail \cite{FM} that every regular matroid is balanced.  In fact,
even more is  true.  Enhancing Feder and Mihail's proof, Choe
\cite{C1,C2} has recently shown the following.
\begin{THM}[Choe \cite{C1,C2}]
Let $\M$ be a sixth--root of unity matroid, and let $e,f\in E(\M)$ be
distinct.  There are sixth--roots of unity $C_{ef}(S)$ for each
$S\subset E$ such that both $S\cup\{e\}$ and $S\cup\{f\}$ are
bases of $\M$, such that
$$\Delta M\{e,f\}(\y)=\left(\sum_{S}C_{ef}(S)\y^{S}\right)
\left(\sum_{S}\overline{C_{ef}(S)}\y^{S}\right).$$
\end{THM}
Since the factors on the right--hand side are complex conjugates when
all the $y_{e}$ are real,  Theorem 5.2 shows that for a sixth--root
of unity matroid $\M$ and distinct $e,f\in E(\M)$, the Rayleigh
difference $\Delta M\{e,f\}(\y)$ is nonnegative for \textbf{any real
values} of the variables $\y$ -- positive, negative, or zero.  We
shall call such matroids  \emph{strongly Rayleigh}.  
\begin{PROP}
Let $\M$ be a strongly Rayleigh matroid on the set $E$.
Then, for all distinct $e,f,g\in E$ and $\y\in\RR^{E}$,
$$|\Theta M\{e,f|g\}|\leq
2\sqrt{\Delta M_{g}\{e,f\}\Delta M^{g}\{e,f\}}.$$
\end{PROP}
\begin{proof}
For a strongly Rayleigh matroid $\M$ and real numbers $\y\in\RR^{E}$
we have $\Delta M\{e,f\}\geq 0$.  Considered as a quadratic polynomial
in $y_{g}$, this does not change sign for $y_{g}\in\RR$, and therefore
it has a nonpositive discriminant.  This gives the stated inequality.
\end{proof}

Arguments directly analogous to those in Section 3 suffice to
prove the following, and the details are omitted.
\begin{PROP}
The class of strongly Rayleigh matroids is closed by taking
duals, minors, and $2$--sums.
\end{PROP}

\begin{THM}
A binary matroid is strongly Rayleigh if and only if it is regular.
\end{THM}
\begin{proof}
It is a theorem of Tutte that a binary matroid is regular if and
only if it does not contain $\F_{7}$ or $\F_{7}^{*}$ as a minor
(Theorems 13.1.1 and 13.1.2 of Oxley \cite{Ox}, for example).
Regular matroids are strongly Rayleigh by Theorem $5.2$.
By Proposition 5.4, to prove the converse it suffices to show that
$\F_{7}$ is not strongly Rayleigh.  Label the elements of $E(\F_{7})$ 
by $\{1,\ldots,7\}$ corresponding to the columns of the
representing matrix
$$\left[\begin{array}{lllllll}
1 & 0 & 0 & 0 & 1 & 1 & 1\\
0 & 1 & 0 & 1 & 0 & 1 & 1\\
0 & 0 & 1 & 1 & 1 & 0 & 1\end{array}\right]$$
over $GF(2)$.  To simplify notation we will write $F_{1}^{26}$
instead of $(F_{7})_{1}^{2,6}$, \emph{et cetera}.  With the 
substitutions $y_{3}=y_{5}=2$ and $y_{4}=y_{7}=-1$ and $y_{6}=t$,
we have $F_{126}=0$, $F_{12}^{6}=F_{16}^{2}=F_{26}^{1}=2$, $F_{1}^{26}
=-8$, $F_{2}^{16}=F_{6}^{12}=1$, and $F^{126}=-4$.  Therefore
\showon
\Delta F\{1,2\}
&=& F_{1}^{2}F_{2}^{1}-F_{12}F^{12}\\
&=& (2t-8)(2t+1)-(2)(t-4)=4t(t-4).
\showoff
For any $0<t<4$ we have $\Delta F\{1,2\}<0$, so that
$\F_{7}$ is not strongly Rayleigh.
\end{proof}

In the case of graphs, Theorem 5.2 specializes to the following
combinatorial identity:\ see also equation (2.34) of Brooks, Smith, Stone,
and  Tutte \cite{BSST}, Theorem 2.1 of Feder and Mihail \cite{FM}, and
several of the identities in Section 3.8 of Balabanian and Bickart \cite{BB}.
\begin{THM}
Let $G=(V,E)$ be a connected (multi)graph, and let $\G$ be the graphic
matroid of $G$.  For distinct $e,f\in E$, fix arbitrary orientations of
$e$ and $f$, and for each $S\subset E$ such that both $S\cup\{e\}$
and $S\cup\{f\}$ are spanning trees of $G$, let $C_{ef}(S):=\pm 1$
according to whether or not $e$ and $f$ are directed consistently
around the unique cycle of $S\cup\{e\}\cup\{f\}$.  Then
$$G_{e}^{f}(\y)G_{f}^{e}(\y)-G_{ef}(\y)G^{ef}(\y)=
\left(\sum_{S}C_{ef}(S)\y^{S}\right)^{2}.$$
\end{THM}
\noindent
A combinatorial proof of this fact is greatly to be desired.

Chavez \cite{CL} has shown that every finite projective geometry
is negatively correlated.  More generally:
\begin{PROP}
If a matroid admits a $2$--transitive group of automorphisms then
it is negatively correlated.
\end{PROP}
\begin{proof}
Let $\M=(E,\B)$ be a matroid of rank $r$ on $m\geq 2$ elements which has a 
$2$--transitive automorphism group, and let $M=M(\one)$, \emph{et 
cetera}.  Let $e,f\in E$ be distinct.  By transitivity of the automorphism
group, $mM_{e}=mM_{f}=rM$.  By $2$--transitivity of the automorphism group,
$m(m-1)M_{ef}=r(r-1)M$.  Thus
$$\Delta M\{e,f\}=M_{e}M_{f}-M_{ef}M=
\frac{M^{2}r(m-r)}{m^{2}(m-1)}\geq 0$$
since $r\leq m$.
\end{proof}

Aaron Williams has recently computed that the finite projective planes
of orders $3$ and $4$ are balanced (personal communication, June 2003).
In the other direction:
\begin{PROP}
Every finite projective geometry is not a HPP matroid.
\end{PROP}
\begin{proof}
Every finite projective geometry contains a finite projective plane as
a minor, so it suffices to prove that finite projective planes are not
HPP matroids.  In fact, a projective plane of order $q$ fails the
condition RZ[$q+1$], as can be seen by taking $S\subseteq E$ to be a
line of the plane and $\y\equiv\one$.  Then the relevant polynomial is
$Ax^{2}+Bx+C$ with
\showon
A &=& q^{2}\binom{q+1}{2}=\frac{(q+1)q^{3}}{2}\\
B &=& (q+1)\left[\binom{q^{2}}{2}-q\binom{q}{2}\right]=
\frac{(q+1)q^{3}(q-1)}{2}\\
C &=& \binom{q^{2}}{3}-(q+1)q\binom{q}{3}=
\frac{(q+1)q^{3}(q-1)^{2}}{6}
\showoff
which has discriminant $-(q+1)^{2}q^{6}(q-1)^{2}/12$, and thus has
non--real zeros.  Theorem 4.5 thus implies that a projective plane
of order $q$ is not a HPP matroid.
\end{proof}

In Section 10.5 of \cite{COSW}, the question is raised whether or not
every transversal matroid is a HPP matroid.  Numerical experiments 
support this idea for transveral matroids of rank three, but we can
no longer hope for much more than this:
\begin{PROP}
There is a transversal matroid of rank $4$ which is not balanced.
\end{PROP}
\begin{proof}
Let $\L$ be the matroid on the set $E=\{1,2,\ldots,10,e,f\}$
for which the bases are the transversals to the four sets
$\{1,2,3,4,f\}$, $\{5,6,7,f\}$, $\{8,9,10,f\}$, and 
$\{1,2,3,5,6,8,9,e,f\}$.  A direct computation shows that
$L_{e}=80$, $L_{f}=168$, $L_{ef}=33$, and $L=436$, so that
$\Delta L\{e,f\}=-948<0$.
\end{proof}

\begin{PROP}
Every matroid with at most $7$ elements is Rayleigh. 
\end{PROP}
\begin{proof}[Sketch of proof]
Since the Rayleigh property is preserved by duality, it suffices to
consider matroids $\M$ for which $\mathrm{rank}(\M)\leq|E(\M)|/2$.
In Table 2 and Appendix A.2 of \cite{COSW}, nine matroids with $7$
elements and rank $3$ are identified as the only matroids
with $|E|\leq 7$ and rank $\leq 3$ which are not known to be
HPP matroids.  (Five are known not to be HPP, four are of unknown status.)
The other small matroids, being HPP, are Rayleigh by Corollary $4.9$.
One of the nine suspicious matroids is the Fano matroid $\F_{7}$, which
was shown to be Rayleigh in the proof of Theorem 3.8.  For each of the
eight remaining matroids a direct \textsc{Maple}--aided calculation
showed that it is Rayleigh.

For example, take the case of $\P'_{7}$, the rank $3$
matroid on $\{1,2,\ldots,7\}$ with three--point lines $\{1,2,6\}$,
$\{2,3,4\}$, $\{1,3,5\}$, and $\{5,6,7\}$.  One finds that $\Delta 
P'_{7}\{e,f\}(\y)$ is a polynomial with nonnegative coefficients 
except when $\{e,f\}$ is one of $\{1,4\}$, $\{1,7\}$, $\{2,5\}$,
or $\{3,6\}$.  The first and second of these cases are equivalent
by an automorphism of $\P'_{7}$, as are the third and fourth, so we
need only consider $\{1,4\}$ and $\{2,5\}$.  In these two cases one
finds that $\Delta P'_{7}\{e,f\}(\y)$ is a positive sum of
monomials and squares of binomials, similar in form to $\Delta 
A_{8}\{e_{7},e_{8}\}$ calculated in the proof of Theorem 3.8.
Thus, $\P'_{7}$ is Rayleigh.

The seven other relevant matroids are handled analogously, and
all are found to be Rayleigh.
\end{proof}

\begin{THM}
The class of balanced matroids is not closed by taking $2$--sums.
\end{THM}
\begin{proof}
By Corollary 3.7 it suffices to give an example of a matroid which is
balanced but not Rayleigh.  The matroid $\J'$ represented over $\RR$ by
the matrix
$$\left[\begin{array}{llllllll}
1 & 1 & 1 & 1 & 1 & 1 & 1 & 3\\
0 & 1 & 0 & 0 & 2 & 0 & 0 & 1\\
0 & 0 & 1 & 0 & 0 & 2 & 0 & 1\\
0 & 0 & 0 & 1 & 0 & 0 & 2 & 3 \end{array}\right]$$
is such an example.  Let $E(\J')=\{1,\ldots,8\}$ corresponding to the
columns of the above matrix.  By Proposition 5.10, every proper minor of
$\J'$ is Rayleigh, so it suffices to show that $\J'$ is negatively
correlated but not Rayleigh.  Straightforward \textsc{Maple}--aided
calculations show that $\J'$ is negatively correlated:\ the value of
$\Delta J'\{e,f\}(\one)$ is given in the $(e,f)$--th entry of this matrix:
$$\left[\begin{array}{rrrrrrrr}
*   & 100 & 100 & 120 & 100 & 100 & 80  & 0\\
100 & *   & 25  & 50  & 225 & 75  & 50  & 100\\
100 & 25  & *   & 50  & 75  & 225 & 50  & 100\\
120 & 50  & 50  & *   & 50  & 50  & 224 & 80 \\
100 & 225 & 75  & 50  & *   & 25  & 50  & 100\\
100 & 75  & 225 & 50  & 25  & *   & 50  & 100 \\
80  & 50  & 50  & 224 & 50  & 50  & *   & 120 \\
0   & 100 & 100 & 80  & 100 & 100 & 120 & *
\end{array}\right]$$
(the diagonal entries are undefined).
However, if the elements are assigned weights
$y_{2}=y_{3}=y_{4}=t$ and $y_{5}=y_{6}=y_{7}=1$, then
$$\Delta J'\{1,8\}(\y)= (t+1)^{3}(t-1)(t^{2}+t-1)$$ and therefore
$\Delta J'\{1,8\}<0$ if $(\sqrt{5}-1)/2<t<1$.  Therefore, $\J'$ 
is not Rayleigh.
\end{proof}
(The matroid $\J'$ in the proof of Theorem $5.11$ is similar in
structure to the sixth--root of unity matroid called $\J$ by
Oxley \cite{Ox}.)

\section{Open Problems.}

The class of Rayleigh matroids is naturally motivated by generalization
of a physically intuitive property, and it has some useful
structure and relevance to other interesting classes of matroids.
There are still many unsolved problems concerning these ideas, among
them the following.

With regard to finding more examples of Rayleigh matroids:\\
$\bullet$\ Is every matroid of rank three Rayleigh?\\
Or, somewhat less ambitiously:\\
$\bullet$\ Is every finite projective plane a Rayleigh matroid?\\
Theorems 3.8 and $5.11$ and Proposition 5.9 show that we can not
hope for all matroids of rank $4$ to be Rayleigh.\\
$\bullet$\  Characterize the class of rank $4$ Rayleigh matroids
by means of excluded minors.

With Theorem 3.8 in mind:\\
$\bullet$\ Characterize the class of ternary Rayleigh matroids
by means of excluded minors.\\
$\bullet$\  Characterize the class of $GF(4)$--representable
Rayleigh matroids by means of excluded minors.\\
Proposition 4.1 provides a starting point for these problems, from
which the method of proof of Theorem 3.8 could be launched.  Completing
either of these projects will require a substantial amount of work,
but should be well worth it.

Concerning the spectrum of conditions between the HPP and Rayleigh
property:\\
$\bullet$\ Is there a Rayleigh matroid which is not
\textup{LC[$4$]}?

Regarding Theorem 5.5:\\
$\bullet$\ Are there strongly Rayleigh matroids which are
not HPP, or not sixth--root of unity?\\
$\bullet$\ Is every HPP matroid strongly Rayleigh?

Finally, in order to better understand the enumerative combinatorics
of graphs:\\
$\bullet$\ Find a combinatorial (bijective) proof of Theorem 5.6.

\end{document}